%Paper: q-alg/9505011
%From: fronsdal <fronsdal@physics.ucla.edu>
%Date: Mon, 15 May 95 16:10:17 0800

%\\
%Title: Cohomology of Quantum Groups.
%Author: Christian Fronsdal.
%Comments: 14 pages, plain Tech.
%\\
%Abstract:Lecture notes. Introduction to the cohomology of algebras, Lie
%algebras, Lie bialgebras
%and quantum groups. Contains a new derivation of the classification of
%classical r-matrices in terms of
%deformation cohomology, and a calculation of the esoteric forms of quantum
%gl(N) by deformation theory.
%\\

\font\stepthree=pbkd scaled\magstep4
%%%\font\stepthree=cmr10 scaled\magstep4
\magnification=\magstep1
\settabs 18 \columns
%paper, date, \b, \q, \r, \ce, ,ve,
%\YB, \UT,
%\hoffset=1.00truein
%\voffset=1.00truein
\hsize=16truecm
\baselineskip=17 pt

\def\b{\bigskip}
\def\bb{\bigskip\bigskip}

\def\no{\noindent}
\def\r{\rightline}
\def\ce{\centerline}
\def\ve{\vfill\eject}
\def\q{\quad}
\def\YB{Y\hskip-1mmB}

\def\UT{U\hskip-0.8mm T}
\def\eqa{\eqalign}
\def\r{\rightline}
\font\got=eufm8 scaled\magstep1
\def\g{{\got g}}
\def\bow{\triangleright \hskip-0.5mm \triangleleft }
\def\gg{\g$\,\bow\, $\g}
\def\ge{\g \enskip}
\def\ep{\epsilon}
\def\G{\Gamma}
\def\pa{\partial}

\def\harr#1#2{\smash{\mathop{\hbox to .25 in{\rightarrowfill}}
 \limits^{\scriptstyle#1}_{\scriptstyle#2}}}
\def\varr#1#2{\llap{$\scriptstyle #1$}\left\downarrow
 \vcenter to .25in{}\right.\rlap{$\scriptstyle #2$}}
\def\diagram#1{{\normallineskip=8pt
 \normalbaselineskip=0pt \matrix{#1}}}

\def\today{\ifcase\month\or January\or February\or March\or April\or
May\or June\or July\or August\or September\or October\or November\or
December\fi \space\number\day, \number\year }

%\r \today

\parindent=0pt

%Ecriture des corps de nombres

\def\Cit{\hbox{\it l\hskip -5.5pt C\/}}

\def\Crm{\hskip0.5mm \hbox{\rm l\hskip -5.5pt C\/}}
\def\w{\wedge}
\def\ep{\epsilon}
\def\G{\Gamma}
\def\pa{\partial}
\def\D{\Delta}
\def\A{\hbox{$A(\langle q\rangle,a)$}}

{\ce {\stepthree   COHOMOLOGY}}
\b
{\ce{\stepthree AND QUANTUM GROUPS}}

\b
{\ce {C. Fr\o nsdal}}

{\ce {Physics Department, University of California, Los Angeles CA 90024, USA}}
\b
\q\q {\bf Preamble.}
\b
\q{\bf 1. Deformations of Quantum gl(n).}

\q\q A.  MULTIPARAMETER QUANTUM $gl(N)$.

\q\q B. DEFORMATIONS.

\q\q C. THE CLASSICAL LIMIT.
\b

\q{\bf 2. Deformation of twisted, simple Lie bialgebras.}

\q\q A. THE PROBLEM.

\q\q B. CALCULATION OF $H^2($\g$^*, \Cit$).

\q\q C. THE HIGHER ORDERS.

\b

\q{\bf 3. Cohomology.}

\q\q A. LIE ALGEBRAS.

\q\q B. LIE BIALGEBRAS.

\q\q C. MANIN ALGEBRA.

\b

\q{\bf 4. Quantization.}

\q\q A. SYMPLECTIC STRUCTURE AND QUANTIZATION.

\q\q B. *-PRODUCTS AND ABSTRACT ALGEBRAS.

\q\q C. COMPATIBILITY.

\q\q D. FURTHER COMMENTS ABOUT QUANTIZATION.

\q\q E. THE OPERATORS $d$ AND $\partial$ ON $U_{\hbar}$.

\q\q F. THE R-MATRIX.
\b
\q{\bf 5. Deformations of twisted quantum groups.}
\b
\q\q {\bf References.}
\ve

\ce {\bf Preamble.}

\b
\def\A{\hbox{$A(\langle q\rangle,a)$}}

\q Section 1 is an abreviated version of a preprint
with the same title by C. Fr\o nsdal and A. Galindo [17].
A  new and unpublished result by the same authors can be reported here.
The representations of the braid group
that are obtained from R-matrices associated with
multiparameter quantum $gl(N)$ are factored through the Hecke algebra.
In other words, the spectrum of the generator consists
of only two points; we refer to this as the Hecke condition.
It turns out that the deformations that preserve the
braid relation (Yang-Baxter) automatically preserve
the Hecke condition. The feature of the R-matrix that naturally
associates it to $gl(n)$ is thus preserved by deformations.
Section 2 \break is an alternative determination of the classical r-matrices
for simple Lie algebras first found by Belavin and Drinfeld.
The methodology consists of calculating the deformations of the simplest
class of coboundary Lie bialgebras (those that respect the Cartan
subalgebra and are
related to the twisted quantum groups).
They are mostly rigid to essential deformations, but on special
surfaces in the space of parameters one finds essential deformations
that turn out to reproduce the entire panorama of Belavin-Drinfeld  \break
r-matrices. Section 3 reviews the tools used,  cohomology
on Lie algebras and the double complex on Lie bialgebras.
There follows an application to the Manin triplet; we prove that,
as a Lie algebra, \gg$^*$ \enskip is isomorphic to \g $\,\,\oplus$ \g \enskip
in all cases. Section 4 deals with quantization and returns to the problem
of deformations of quantum groups in the general case. A lecture
on the Universal T-matrix, in which a solution was offered to the problem
of exponentiation on quantum groups, is not covered by these notes. The
material
may be found in ref.10.
\bb

\ve
\ce{{\bf 1. Deformations of Quantum gl(n).} }
\b
A.  MULTIPARAMETER QUANTUM $gl(N)$.%1A

\q  Belavin and Drinfeld [1] classified the r-matrices (structures of
coboundary
Lie bialgebra)
associated to simple Lie algebras, both finite and affine.
A program of ``quantization" of Lie algebras
proposed by Drinfeld [2] would promote the classical structures to bialgebra
deformations of enveloping algebras.  Standard forms of such
``quantum deformations" of the simple affine
Lie algebras were obtained by Jimbo [3].
Here we are interested only in the finite case; that is, constant r- and
R-matrices.  So far there is no general classification of quantized
Lie algebras.

\q Deformation theory applied at the classical point [4] is difficult, since
the obstructions appear only in the second order.  But a large family of exact
quantum deformations of $gl(N)$ is known [5], with $1 + N(N-1)/2$ parameters.
These algebras are rigid (with respect to essential deformations)
 at generic points in parameter space, even
to first order deformations, but essential deformations exist on
algebraic surfaces of lower dimension, for $N>2$.
 The determination of all the
first order deformations, presented here, goes far towards a
complete classification of all formal and/or exact deformations.  All first
order deformations are combinations of ``elementary" deformations, and all the
elementary deformations are exact. [A large class of exact
deformations were described in ref.9; I learned at this School that a
one-parameter family belonging to
the same class had already been discovered by Cremmer and Gervais, ref.18.]

\q We consider the free associative algebra ${\cal F}_x$ generated by $(x^i),
\,\, i =
\cdots, N$, and the ideal ${\cal F}_{x0}$ generated by
$$x^ix^j - q^{ij}x^jx^i, \,\, i,j=1,\cdots ,N ,\eqno(1.1)
$$
\no in which the $q$'s are taken from a field K with characteristic 0, with
$$
q^{ij}q^{ji} = 1, \,\, q^{ii} = 1, \,\, i,j = 1, \cdots ,N . \eqno(1.2)
$$
\no We call {\bf quantum plane} the associative algebra
${\cal F}_x/ {\cal   F}_{x0}$;
that is, the associative algebra generated by the $x$'s with relations
$$
x^ix^j - q^{ij}x^jx^i = 0, \,\, i,j = 1, \cdots ,N . \eqno(1.3)
$$
\no Similarly, the quantum anti-plane
${\cal F}_{\theta}/ {\cal F}_{\theta 0}$ is
an associative algebra generated by N elements
$(\theta^i), \,\, i = 1, \cdots
, N$, with relations
$$
\theta^i \theta^j + r^{ij} \theta^j \theta^i = 0, \,\, i,j = 1, \cdots ,N ,
\eqno(1.4)
$$
\no in which the $r$'s are parameters from $K$
satisfying the same relations as
the $q$'s, Eq. (1.2).  Let $V$ denote the linear
vector space over $K$ spanned by
the $x$'s (or by the $\theta$'s).
\b
\no {\bf Definition}.  A {\bf generalized symmetry} is an element $P$ of
${\rm End}(V \otimes V)$ that satisfies the Hecke condition
$$
(P-1)(P+a) = 0, \eqno(1.5)
$$
\no for some $a \in K, a \not= -1,0$.  We shall say that the tensor $xx =
(x^ix^
j)$ is
$P$-symmetric, and that the tensor $\theta \theta = (\theta^i \theta^j)$ is
$P$-antisymmetric, if
$$xx(P-1) = 0, \,\,\, \theta \theta(P+a) = 0 . \eqno(\hbox{1.3'-4'})$$
\b
\q Let $P_{12}$ be the operator on $V \otimes V \otimes V$ that acts as $P$ on
the two first factors; the braid relation is
$$
P_{12} P_{23} P_{12} = P_{23} P_{12} P_{23} . \eqno (1.6)
$$
\b
\no{\bf Theorem 1}. Given relations (1.3) and (1.4), with parameters $q$ and
$r$ subject to the conditions (1.2),
the following two statements are equivalent:
$${\rm (a)} \,\,\,
   r^{ij} = aq^{ij}, \,\, i < j, \,\, i,j = 1, \cdots, N ; \eqno(1.7)
$$
\no{\rm (b)} There exists $P$ in ${\rm End}(V \otimes V)$
satisfying (1.5-6) ,
such that the
relations (1.3-4) coincide respectively with Eqs. (1.3'-4');  it is unique
up to a permutation of the basis.
\b
\q Let $P$ be a generalized symmetry of dimension $N$.
Consider the algebra ${\cal F}_x$  freely generated
by $(x^i), \,\, i=1, \cdots ,N$, with the ideal
${\cal F}_{x0}$ generated by

$$
(xx(P-1))^{ij}, \,\, i,j = 1, \cdots, N;\eqno(1.8)
$$
\no and the algebra ${\cal F}_{\theta}$ generated by
$(\theta^i) \,\, i= 1, \cdots , N$,
with the ideal ${\cal F}_{\theta 0}$ generated by

$$
\theta \theta (P+a).\eqno(1.9)
$$
\no Let ${\cal F}$ be the algebra generated by
$(x^i)$ and $(\theta^i), \,\, i= 1, \cdots , N$ with relations
$$
ax \theta = \theta xP . \eqno(1.10)
$$
\no This algebra contains ${\cal F}_x$ and ${\cal F}_{\theta}$ as subalgebras
and the ideals ${\cal F}_{x0}$ and ${\cal F}_{\theta 0}$ are thus canonically
identified with subsets of ${\cal F}$.
\b
\no {\bf Theorem 2.}  Suppose that the generalized symmetry $P$ satisfies the
braid relation.  Let $X$ denote the linear span of the $x$'s and $\Theta$ the
linear span of the $\theta$'s, then the statements
$$
\eqalignno{
{\cal F}_{x0} \Theta = \Theta {\cal F}_{x0}&&(1.11)\cr
{\cal F}_{\theta 0} X = X {\cal F}_{\theta 0}&&(1.12)\cr}
$$
\no hold in ${\cal F}$.  Conversely, if
both (1.11) and (1.12) hold, then $P$
satisfies the braid relation.
\b
\no {\bf Proof.}  Eqs. (1.11), (1.12) are equivalent, respectively, to
$$
({\rm braid})_{123} (P_{12}-1) = 0, \quad
({\rm braid})_{123} (P_{12}+a) = 0, \eqno(1.13)$$
with
$$
({\rm braid})_{123}:= P_{12} P_{23} P_{12} - P_{23} P_{12} P_{23}.
\eqno(1.14)
$$
\b
\no {\bf Definition.}  Let $\langle q\rangle$ stand
for a set of parameters $(q^{ij}), \,\,
i,j = 1, \cdots, N$, satisfying $q^{ij}q^{ji} = 1$ and $q^{ii} = 1$; and $a$ an
additional parameter, all in the field $K$.
Let $r^{ij} = aq^{ij}$ for $i<j, \,\,
r^{ii} = 1$ and $r^{ji} = q^{ji}/a$ for $i<j$.  The {\bf standard quantum
algebra} {\A}  is generated by the $x$'s and the $\theta$'s, with
relations (1.3-4) and (1.10).  More generally,
for any generalized symmetry $P$,
the algebra $A(P)$ generated by the $x$'s and the $\theta$'s, with relations
(1.3'-4') and (1.10), will be called a {\bf quantum P-algebra}
if the conditions (1.11) and (1.12) hold.
\ve
\q The purpose of this paper is to study deformations of {\A}  in the
category of quantum P-algebras.  The quantum pseudogroup (in the sense of
Woronowicz [6]) associated to $P$ is the unital algebra generated by the matrix
elements of an $N$-by-$N$ matrix $T$, with relations
$$
[P,T\otimes T]=0.\eqno(1.15)
$$
\no It is the algebra of linear automorphisms of $A(P)$; that is, the
set of mappings
$$
(x,\theta,T)\rightarrow x\otimes T, \,\,\, \theta \otimes T\eqno(1.16)
$$
\no that preserve the relations (1.3'-4', 1.10) of $A(P)$.
It is related, via duality,
to a quantum group in the sense of Drinfeld [2].  Twisted, quantum $gl(N)$ [5]
corresponds to ${\rm Aut}\A$.  The deformations of this quantum group are
 in 1-to-1 correspondence with the deformations of the standard quantum
algebra \A.

\q The next section defines the deformations that will be calculated in this
paper.
 It turns out that {\A} is rigid for parameters in general position.
Interesting
nontrivial deformations (even exact ones) do exist on certain algebraic
surfaces in parameter space, for $N>2$.    The existence of an unexpected
special case
that requires $a^3 = 1$ deserves some attention.

\bb
B. DEFORMATIONS.%1B

\q Henceforth,  $P$ will denote the generalized symmetry associated with \A.
A {\bf formal deformation} of {\A}   is here a quantum $P(\epsilon)$-algebra
with $P(\epsilon)$ a formal power series in an indeterminate $\epsilon$
$$
P(\epsilon)=P + \epsilon P_1 + \epsilon^2 P_2 + \cdots ,\eqno(1.17 )
$$
\no that satisfies the Hecke condition with the parameter $a$ independent of
$\epsilon$, and such that (1.11-12) hold.  In this case we shall say that
$P(\epsilon)$ is a formal deformation of $P$.  A deformation is {\bf exact} if
the series $P(\epsilon)$ has a nonvanishing radius of convergence.
 If $P(\epsilon)$ is a formal deformation of $P$, then
$$
P(\epsilon,1)=P+\epsilon P_1\eqno(1.18)
$$
\no is a first order deformation.  More generally, a {\bf first order
deformation} is defined as a formal deformation except that one sets
$\epsilon^2 = 0$.  A first order deformation is not necessarily the first two
terms of a formal deformation.  For example, at the classical point, where $a$
and all the $q$'s are equal to unity, the braid relation is moot in first
order.  (Recall that the condition that defines a formal deformation is in this
case the classical Yang-Baxter relation, which is second order.)  For this
reason, the concept of a first order deformation is of no use at the classical
point.  In contrast with this, we shall find that, at general position in
parameter space, {\A}   is rigid with respect to first order deformations,
which implies rigidity under formal deformations.  The case $a=1$, with the
$q$'s in general position, is in this respect intermediate; it is best treated
separately.

Two types of deformations (and combinations of them) will be considered
trivial.  A linear transformation, with coefficients in $K[\epsilon]$,
$$
x^i \rightarrow x^i + \epsilon x^j A_j^i + \cdots , \,\, \theta^i \rightarrow
\theta^i + \epsilon \theta^j A_j^i + \cdots , \eqno(1.19)
$$
\no induces a trivial, formal deformation of $P$.  A variation of the $q$'s
$$
q^{ij} \rightarrow q^{ij} + \epsilon \delta q^{ij} + \cdots \eqno(1.20)
$$
\no will also be considered trivial.  A deformation that is not trivial is
called essential.  We shall classify the equivalence classes, with respect to
the transformations (1.19-20), of first order deformations.

\b
\no {\bf Theorem 3.}  If $a \not= 1$, then each equivalence class of first
order deformations contains a unique representative with the property that
$(P_1)_{kl}^{ij} = 0$ for every index set $i,j,k,l$ that contains no more
than two different numbers.
\b
\q
The first order deformation of $P$ induced by (1.19) is
$$
P_1 = PZ - ZP, \quad Z:= A \otimes 1 + 1 \otimes A,\eqno(1.21)
$$
\no or more explicitly
$$
(P_1)_{kl}^{ij} = a(\hat{q}^{lk} - \hat{q}^{ij}) Z_{kl}^{ji} + (1-a)[(k<l)
 - (i<j)] Z_{kl}^{ij}, \eqno(1.22)
$$
\no where
$$
{
\hat{q}^{ij}:=\cases{
q^{ij}&{\rm if} $i<j$,\cr
q^{ij}/a&{\rm if} $i\geq j$,\cr}
\qquad
(i<j):=\cases{
1&{\rm if} $i<j$,\cr
0&{\rm
otherwise}.\cr}
}
\eqno(1.23)
$$

Preservation of the Hecke condition (1.5) under first order deformations is
equivalent to requiring that
$$
\eqalign{q^{ij} (P_1)_{lk}^{ji} + \hat{q}^{kl} (P_1)_{kl}^{ij} &=0 \cr
(1-a) q^{ij} (P_1)_{lk}^{ji} =(a-1) \hat{q}^{kl} (P_1)_{kl}^{ij}
&= (P_1)_{lk}^{ij} + aq^{ij} \hat{q}^{kl} (P_1)_{kl}^{ji},
 \,\, i \leq j, \,\,
k \leq l.\cr}\eqno(1.24)
$$
\no The main difficulty is to extract the conditions on $P_1$ imposed by the
braid relation.  We found that the strategy made available by Theorem 2
simplifies this task.  Both conditions, (1.11) and (1.12), must be invoked.  We
leave out the details.
\b

\no {\bf Definition.}  An elementary first order deformation is one in which
some $(P_1)_{kl}^{ij}$ is non-zero for just one unordered pair $i,j$
and just one unordered pair $k,l$.
\b
\no {\bf Theorem 4.}  Suppose $a \neq \pm 1, 0 $.
There are two series of elementary,
first order, essential deformations. The ``principal series"
is described first. Let $i,j$ be any index pair with, either (case 1)
$ k + 1 = i \leq j = l - 1$, or else (case 2) $i+1 = k \leq l = j-1$.  Let
$P_1=0$, except that
$$
(P_1)_{lk}^{ij} = -a\hat q^{ij}q^{kl} (P_1)_{kl}^{ji} \not= 0.\eqno(1.25)
$$
\no This defines an elementary deformation if and only if the parameters
satisfy the conditions
$$
q^{im} q^{jm} q^{mk} q^{ml} = \cases{a^x,&case 1, \cr
                                    a^{-x}, &case 2; \cr}\,\, \,
x = \delta_m^i - \delta_m^j, \, m = 1,2 \ldots N. \eqno(1.26)
$$
\no The ``exceptional series" of elementary, first order deformations exists
only if $a^3 = 1$.
Let $i,j,k$ be neighbors in the natural numbers, with $i + 1 = j$. Let $P_1=0$
, except that either
$(P_1)_{kk}^{ij} = -aq^{ij}(P_1)_{kk}^{ji} \not= 0$, or else
$(P_1)_{ij}^{kk} = -q^{ij}(P_1)_{ji}^{kk} \not= 0$, but not both.
This defines an elementary deformation if and only if the parameters satisfy
$$
\eqalign {&(P_1)_{kk}^{ij} \not= 0:\,\,\,
         (q^{km})^2q^{mj}q^{mi} = a^x,\,\,\, x = \delta_{mi} - \delta_{mj};
           \cr
 &(P_1)_{ij}^{kk} \not= 0: \,\,\,
         (q^{km})^2q^{mj}q^{mi} = a^x,\,\,\,
x = \pm (\delta_{mk} - \delta_{mi}). \cr} \eqno(1.27)
$$
\no The two signs in the last line apply when $k = i-1, k = j+1$, respectively.
There are no other first order, elementary deformations.  The elementary
deformations are formal and exact.

\bb
C. THE CLASSICAL LIMIT.%1D

\q All the results obtained
here,
for twisted quantum $gl(N)$,
have direct
application to twisted, quantum $sl(N)$.
The connection between these two was
explained by Schirrmacher in [5] and is discussed also in [10].
In this section
we shall take the classical limit and confront our results for $gl(N)$
with the
classification, by Belavin and Drinfeld [1], of the classical r-matrices
for $sl(N)$. Strictly, this is possible only under additional conditions
on the
parameters, namely
$$
\prod_iq^{ij}a^j = a^{(N+1)/2};\eqno(1.28)
$$
\no one may therefore assume that these relations hold,
although they do not interfere
directly with the following calculations.

\q The deformed quantum P-algebras of the principal series are semiclassical.
The classical r-matrix is defined by expanding the parameters,
$$
a = 1 + h, \, q^{ij} = 1 + hp^{ij}, \, i<j, \eqno(1.29)
$$
\no and the R-matrix,
$$
R_{lk}^{ij}:=P_{kl}^{ij} + \epsilon (P_1)_{kl}^{ij}, \eqno(1.30)
$$
\no in powers of $h$,
$$
R = 1 - hr_\epsilon + O(h^2), \, r_\epsilon = r + \epsilon \delta r.
  \eqno(1.31)
$$
\no Here $r$ is the r-matrix for twisted (= multiparameter) $gl(N)$,
$$
r = \sum_{i<j} M_j^i \otimes M_i^j + r_0, \eqno(1.32)$$
\noindent with
$$r_0 := \sum_{i<j} \bigl(p^{ij} M_j^j \otimes M_i^i - (1 + p^{ij})
  M_i^i \otimes M_j^j \bigr), \eqno(1.35)
$$
\no and the perturbation is
$$
h\delta r_{kl}^{ij} = (P_1)_{lk}^{ij}, \,\,\, {\rm or} \, \,\,
  h\delta r = \sum_{i,j,k,l} (P_1)_{lk}^{ij} \,\,M_k^i \otimes
M_l^j.\eqno(1.34)
$$
\no Here $M_k^j$ is the matrix with the unit in row-$k$, column-$j$, all the
rest zero.

 \q We examine the classical limit of an elementary, first-order deformation.
Fix
the notation as in Theorem 4, the expression for $\delta r$ is, up to a
constant,
$$
\delta r = M_k^i \otimes M_l^j - M_l^j \otimes M_k^i .\eqno(1.35)
$$
\no The diagonal matrices $(M_i^i), \, i=1,\ldots,N$, will be taken as a basis
for a ``Cartan subalgebra of $gl(N)$".  The upper triangular matrices form the
subspace of positive roots and the matrices $M_i^j$ with $i-j= \pm 1$ are the
simple roots.  [We have abused the notation by extending the notion of roots
from $sl(n)$ to $gl(n)$ and by introducing both positive and negative ``simple"
roots.]  The conditions on the indices that are spelled out in Theorem 4 insure
that all the roots appearing in (1.35) are simple, and that a positive root is
paired with a negative root and {\it vice versa}.  Let us find out what is the
meaning of the restriction (1.26) on the parameters.

\q Let $\Gamma_1(\Gamma_2)$ be the root space spanned by $M_k^i(M_j^l)$ and
$\tau: \Gamma_1 \rightarrow \Gamma_2$ the mapping defined by $\tau(M_k^i) =
M_j^l$.  Consider the equation
$$
(\tau \alpha \otimes 1 + 1 \otimes \alpha) r_0 = 0, \eqno(1.36)
$$
\no where $\alpha = M_k^i$ and $(\alpha \otimes 1) H \otimes H' =
\alpha(H)H'$ and $ (1 \otimes \alpha)H \otimes H' = \alpha(H')H$,
$H$ and $H'$ in the diagonal subalgebra of $gl(N)$.  We find that this equation
is equivalent to
$$
p^{lm} + p^{km} + p^{mi} + p^{mj} =  \delta_m^j - \delta_m^i,\eqno(1.37)
$$
\no  This
equation is precisely the first order analog of Eq. (1.26), while (1.36) is a
condition (Eq.~(6.7))
of Belavin and Drinfeld [1].
 We have thus established that the conditions (1.26)
on the parameters lift the invariance condition
of ref. [1] to the quantum algebra.
\bb
\ve

\ce{ {\bf 2. Deformation of twisted, simple Lie bialgebras.}}
\b
A. THE PROBLEM.%2A

\q The (constant) classical Yang-Baxter matrices for simple Lie algebras
were already
classified by Belavin and Drinfeld [1].
Here we shall approach the same problem by
the methods of deformation theory; we shall see that the calculations are
straightforward and that this method throws some light on the
classification. The results show a clear relationship to
the deformations of twisted quantum,
$gl(n)$  and suggest that
the deformations of all simple quantum groups can be obtained with
less effort using cohomological methods.

\q This lecture merely shows the calculations; therefore
everything is referred to a basis.
The theoretical background is explained in Section 3.

\q   Let \ge be a simple Lie bialgebra over  \Crm,
 $(L_i)$ a basis, structure tensors $\ep$ and $f$:
$$
[L_i,L_j] = \ep_{ij} = \ep _{ij}^kL_k,
\quad \Delta(L_i) = f_i = f_i^{jk}L_j\otimes L_K.
$$
Here $\ep$ is a evidently a two-form valued in \g,
satisfying the Jacobi identity
$$
(d\ep)_{ijk} = \sum_{(ijk)}\ep_{im}^n\ep_{jk}^m = 0,
$$
and $f$ is a oneform valued in \g$\,\otimes\,$\g.
The cohomology operator of \ge is denoted $d$.
The ``compatibility condition" between the two structures reads
$$
df = 0, \q
(df)_{ij} = L_if_j - L_jf_i - \ep_{ij}^k f_k.
\eqno(2.1)$$
In the case of a coboundary Lie algebra $f$ is exact,
$$
f = dr, \q f_i = [L_i,r], \q r = r^{ij}L_i \otimes L_j.
$$
\q From now on we set $f = dr$. Let \g$^*$ denote the vector space
dual of \g,  dual basis $\Gamma^i$; it becomes a Lie bialgebra with
structure tensors $f$ and $\ep$,
$$
\{\Gamma^i,\Gamma^j\} = f^{ij} = f_k^{ij}\,\G^k, \q
\Delta(\Gamma^i) = \ep^i = \ep_{jk}^i \,\Gamma^j \otimes \Gamma^k.
$$
The cohomoloy operator on the Lie algebra \g$^*$   is denoted $\partial $.
 (Section 3.)
The Jacobi identity for the Lie  bracket $\{,\}$ of \g$^*$ is
$$
\partial f = 0, \q (\partial f)_m =
 f_m^{in}\,f_n^{jk} \,L_i \w L_j \w L_k.
$$
\q Since $f = dr$, and since $d$ and $\partial$ anticommute, it may be
expressed
in terms of the classical Yang-Baxter or Schouten bracket,
$$
\pa f = \pa d r = -d\pa r = 0, \q
\pa r =  f_n^{ij}\,r^{kn}\,L_i \w L_j \w L_k,\eqno(2.2)
$$
or
$$
d(\YB) = 0, \q \YB = \pa r =
[r_{12},r_{13} + r_{23}] + [r_{13},r_{23}].
$$
At the cost of including an invariant, symmetric piece in $r$ we can take $r$
to be $\pa$-closed, $\pa r = 0$. Thus, following Belavin and Drinfeld, we pose
$$
r + r^t = K, \q \pa r = 0,\eqno(2.3)
$$
where $K$ is the Killing form of \g.

\q We can deform in the category of
Lie bialgebras or in the category of boundary Lie
bialgebras.
In the first case we deform $f$, in the second case $r$.
Since the symbol $\epsilon$ is used for the structure tensor
 we use $\hbar$ for the deformation parameter.

\q To deform in the category of Lie bialgebras;
one would  set
$$
f(\epsilon ) = dr + \hbar \, f_1 + ...\,,
$$
and impose associativity and closure,
$$
\partial f(\hbar)  = 0, \q df(\hbar) = 0.
$$
To first order in $\hbar $ ( with $a,b,c \in \,$\g$^*$)
$$
\partial f_1(a,b,c)
= \sum_{cyclic}\left(f_1(a,f(b,c)) + f(a,f_1(b,c)\right)
= 0, \q df_1 = 0.
$$
The relevant cohomologies are $H^1($\g, \g$\w$\g), which is $0$,
 and $H^2($\g$^*,\,\,$\Cit$)$, which is not.

\q Deforming in the category of coboundary Lie algebras;
we set
$$
f(\hbar ) = dr(\hbar),\q
r(\hbar ) = r + \hbar r_1 + ... \,,\eqno(2.4)
$$
and impose associativity only.
To first order in $\hbar$,
$$
r_1 + r_1^t = 0, \q\partial  r_1 = 0,
$$
or more explicitly
$$
\partial r_1 =
 \sum f_m^{ij}\,r_1^{km}\,L_i \wedge L_j \wedge L_k = 0.
\eqno(2.5)
$$

\q  The cohomology
is now $H^2($\g$^*)$  with coefficients in the field,
since $r_1(a,b) \in   K$.
The problem is to calculate $H^2($\g$^*,\Cit)$; that is,
the equivalence classes of essential deformations,
in the category of coboundary bialgebras, of the simplest
Lie bialgebras (that correspond to the twisted quantum groups).

\bb

B. CALCULATION OF $H^2($\g$^*,\Cit)$.%2B

\q  Notation. Choose a Cartan subalgebra {\got h},
a set of positive roots, and a Weyl basis
$(e_{\alpha}, e_{-\alpha})_{\alpha > 0}, (h_i)$.
The structure of \ge is
$$
[h_i,e_{\alpha}] = r_i(\alpha)e_{\alpha},\q
[e_{\alpha},e_{\beta}] = e_{\alpha + \beta}, \,\alpha < \beta, \q
[e_{\alpha},e_{-\alpha}] = r^i(\alpha)h_i.\eqno(2.6)
$$
To make sense of this we must introduce an ordering on the root lattice;
we assume that $\alpha < \beta$ implies
that $-\alpha < -\beta $.
\b
\q For $r$ we have
$$ r = \sum r_0^{ij}\, h_i \otimes h_j +
\sum_{\alpha > 0} e_{\alpha }\otimes e_{-\alpha},\eqno(2.7)
$$
and thus
$$
r + r^t = \sum (r_0^{ij} + r_0^{ji})\,\, h_i \otimes h_j +
\sum_{\alpha  } e_{\alpha }\otimes e_{-\alpha} = K^{ij}\,L_i\otimes L_j.\eqno
(2.8)
$$
Here $r^t$ is the transposed of $r$
and the last expression is the invariant
in \ge $\otimes $ \g. The coefficients
define the Killing form $K$, and the restriction $r_0$
of $r$ to {\got h} is
$$
r_0 = \hat r_0 + (1/2) K_0,\eqno(2.9)
$$
with $\hat r_0 $ antisymmetric and  $K_0$
the restriction of $K$ to {\got h}.
Note that  $r^j(\beta) = K_0^{ij}r_i(\beta)$.

\b
\q  We calculate the coproduct,
$$
\Delta(L_i) = [L_i,r] = \sum f_i^{jk} L_j \wedge L_k,\eqno(2.10)
$$
and see that $
\Delta(h_i) = 0$
while
$$
\eqalign{
\Delta (e_{\beta })& = [e_{\beta},r] =
\sum f_{\beta}^{ij}\,L_i \wedge L_j\cr
&= \sum r_0^{ij}\left([e_{\beta},h_i] \otimes  h_j
+ h_i \otimes  [e_{\beta},h_j]\right)
  + \sum_{\alpha > 0}\left([e_{\beta},e_{\alpha}]
\otimes  e_{-\alpha} + e_{\alpha} \otimes
[e_{\beta}, e_{-\alpha}]\right)\cr
= &\,\,-\hat r_0^{ij}r_i(\beta))e_{\beta} \wedge h_j
+ (1/2)\sum_{\alpha > 0}\left([e_{\beta},e_{\alpha}]
\wedge e_{-\alpha} + e_{\alpha} \wedge
[e_{\beta},e_{-\alpha}]\right)\cr
= &\q \sum \hat r_0^{ij}r_i(\beta ) h_j \wedge e_{\beta}\cr & +
(1/2)(\beta < 0)r^j(\beta)h_j \wedge e_{\beta}
  + (1/2)(\beta > 0)r^j(\beta)e_{\beta} \wedge h_j + ...,\cr
}
$$

\q The coefficients give us $f$ and thus the Lie structure of \g$^*$.
The dual basis will be denoted
$(x^{\alpha } ,y^{-\alpha}, w^i)$.
The structure of \g$^*$ is
$$
\eqalign{
&\{w^i,x^{\alpha}\} = w^i(x^{\alpha})\,x^{\alpha},\,\,
\{w^i,y^{-\alpha}\} = w^i(y^{-\alpha})\,y^{-\alpha},\cr
&\{x^{\alpha},y^{-\beta}\} = 0,\,\,
\{x^{\alpha},x^{\beta}\} = x^{\alpha + \beta},\,\,
\{y^{-\alpha},y^{-\beta}\} = y^{-\alpha - \beta} ,\,\,  \alpha < \beta,\cr}
\eqno(2.11)
$$
The ``weights" $w^i(x^{\alpha }),...$ are in the field;
$$
w^j(x_{\beta}) =  \sum_i\left(\hat r_0^{ij} \mp
(1/2)K_0^{ij}\right)r_i(\beta),\eqno(2.12)
$$
for $\beta $ positive/negative.
\b
\q  Next, the differentials. The zero-forms are closed
and never exact. Consequently, one-forms are never exact,
$B^1 = 0, H^1 =  Z^1 \subset C^1$. If $A$ is a one-form,
$$
A = A^ih_i + A^{\alpha}e_{\alpha},
$$
then
$$
\partial A = f_k^{ij} A^k L_i \wedge L_j =
2\sum w^i(\alpha)A^i h_i \wedge e_{\alpha}
+ \sum_{\alpha, \beta > 0}A^{\alpha + \beta}
e_{\alpha} \wedge e_{\beta}.\eqno(2.13)
$$
So $\partial A = 0$ implies that
$$
A = \sum_{simple}A^{\alpha} e_{\alpha}.\eqno(2.14)
$$
The sum is over simple roots $\alpha\,\,$for which $w$ does not vanish,
which means all simple roots.

\q We turn to two-forms,
$$
B = \sum B^{ij} L_i \wedge L_j,
$$
$$\eqalign{
dB& = \sum f_l^{ij}B^{kl} L_i \wedge L_j \wedge L_k\cr
&= 2\sum w^i(\alpha)B^{j,\alpha}
 h_i \wedge e_{\alpha} \wedge h_j
+ 2\sum w^i(\alpha)B^{\beta, \alpha}
h_i \wedge e_{\alpha} \wedge e_{\beta}\cr
&\q \,\,\,\, + \sum_{\alpha \beta > 0} B^{i,\alpha + \beta}
e_{\alpha} \wedge e_{\beta} \wedge h_i
+ \sum_{\alpha \beta > 0} B^{\gamma, \alpha + \beta}
e_{\alpha} \wedge e_{\beta} \wedge e_{\gamma}.
}
\eqno(2.15)$$
Here $\alpha \beta > 0$ means that the two  roots are
either both positive or both  negative.

\q According to (2.13), a two-form is exact if
$$
B^{ij} = 0, \q B^{i\alpha} = w^i(\alpha)A^{\alpha},
\q B^{\alpha \beta} = \biggl \{\eqalign{ &A^{\alpha + \beta},
 \alpha \beta > 0 ,\cr &0 , \q {\rm otherwise}.\cr }
\eqno(2.16)
$$
According to (2.15), a two-form is closed if
$$
B^{i\alpha} = w^i(\alpha)A^{\alpha},\q
B^{\alpha \beta}\left(w^i(\alpha) + w^i( \beta)\right) =
\biggl \{\eqalign{&B^{i,\alpha + \beta},
\,\,\alpha \beta > 0,\,\, \alpha ,\beta \,\,{\rm simple}\cr
 &0\q \q \,\,, \,\,{\rm otherwise}\cr
} \eqno(2.17)
$$

\q The result (2.12) for the weights shows that
$w(\alpha) + w(\beta)$
 cannot vanish unless
$\alpha \beta < 0.$ If $\alpha,\beta $ are both positive,
then $w(\alpha) + w(-\beta) = 0$ is the same as
$$
r^{ij}r_j(\alpha) + r_j(\beta)r^{ji} = 0.\eqno(2.18)
$$

Compare Eq.(1.36). From (2.16) and (2.17) we draw the following conclusion.
\b
{\bf Proposition.} The space of first order essential deformations
of the
bialgebra is
$$
H^2 = Z^2/B^2 = \bigl\{r_1 =
\sum_{\sigma} r_1^{\alpha,\beta} e_{\alpha}
\wedge e_{\beta} + \sum r_1^{ij}h_i \wedge h_j \bigr\},
$$
where $\sigma $ is the set
$$
\sigma = \{(\alpha,\beta);\,\,
\alpha,\beta \,\,{\rm simple},\, w(\alpha) + w(\beta) = 0\};
$$
and $(\alpha,\beta) \in \sigma$ implies that
$\alpha\beta < 0$. The coefficients $r_1^{ij}$
in the second term merely perturb the original parameters in
$r_0$ and these deformations are considered trivial, to be ignored
in the sequel.

\bb
 \ve

C. THE HIGHER ORDERS.%2C

\q
 To order $\hbar^2$ the closure condition is
$$
d\bigl( \YB(r_1) + \partial r_2\bigr) = 0,
$$
or, since $r_1, r_2$ are antisymmetric,
$$
\YB(r_1) + \partial r_2 = 0.
$$
The first term is, according to the Proposition (up to
equivalence),
$$
\eqalign{
\YB(r_1)
= &(1/2)\sum r_1^{\alpha \beta } r_1^{\gamma \delta }
\bigl( [e_{\alpha } ,e_{\gamma }]
 \wedge e_{-\beta } \wedge e_{-\delta }\,
+ [e_{-\beta },e_{-\delta }]
 \wedge e_{\alpha } \wedge e_{\gamma }\cr
&\q  -  [e_{\alpha} ,e_{-\delta }]
 \wedge e_{-\beta } \wedge e_{\gamma }
-  [e_{-\beta },e_{\gamma }]
 \wedge e_{\alpha } \wedge e_{-\delta }
    \bigr) \, .
\cr
}\eqno(2.19)
$$
It is a general result of deformation theory that this quantity is
$\pa$-closed; now we have to find $r_2$ such that $\pa r_2 = -\YB$,
the obstructions to continuing the deformation are thus  in
$H^3($\g$^*,\Cit)$.
Note that the roots that appear in (2.19) are all simple; therefore
$$
[e_{\alpha} ,e_{-\delta }] =
\delta_{\alpha}^{\delta}\,r^i(e_{\alpha} )h_i ,\q
[e_{-\beta },e_{\gamma }] = -\delta_{\beta}^{\gamma}\,
r^i(e_\gamma)h_i.
$$
One part of the YB bracket is thus
$$
\eqalign{
-\sum r_1^{\alpha \beta} &r_1^{\gamma \alpha}r^i(\alpha )
\,h_i  \,\wedge e_{-\beta } \wedge e_{\gamma }
+ \sum r_1^{\alpha \beta} r_1^{\beta\delta}r^i(\beta )
\,h_i \wedge e_{\alpha } \wedge e_{-\delta }\cr
&= 2\sum  r_1^{\alpha \beta }
r_1^{\beta \delta }
r^i(\gamma )
\,h_i\, \wedge e_{\alpha } \wedge e_{-\delta }.\cr
}\eqno(2.20)
$$

To cancel this we need a term $r_{21}$ in $r_2$ that satisfies
$$
r_{21}^{\alpha,-\delta}\bigl(w^i(\alpha) + w^i(-\delta)\bigr) =
2\sum_{\beta} r_1^{\alpha \beta }
r_1^{\beta \delta }
  r^1(\gamma),
\eqno(2.21)
$$
so that, if $w^i(\alpha) + w^i(-\delta)  = 0$, then
there would be an obstruction. But this is impossible;
since $w(\alpha) + w(-\gamma) = 0$  it would imply that
$w(-\gamma) - w(-\delta) = 0$, and thus $ \gamma = \delta$;
but $r^{\gamma\gamma} = 0$. [Remember that, in (2.21),
the indices $\alpha,\beta, \delta$ refer to simple roots.]

\q Turning to the remaining terms in (2.19),
suppose that $\alpha + \gamma $ is a root,
$\alpha + \gamma = \mu $. Then the first term in (2.19)
becomes
$$
\sum r_1^{\alpha \beta}r_1^{\gamma\delta}\,e_{\mu}
\wedge e_{-\beta} \wedge e_{-\delta}.
$$
To cancel it we need, as (2.15) shows, that two of the roots
in this expression add up to a root.
This can only be $-\beta -\delta = - \nu .$
So here is an obstruction;
if $\alpha + \gamma$ is a root and $\beta + \delta$ is not, and
$r_1^{\alpha,\beta}r_1^{\gamma,\delta} \neq 0$,
then there is no remedy.
We must avoid this obstruction by imposing an additional
condition on $r_1$: in the formula
$$
r_1 =
\sum_{\sigma} r_1^{\alpha,\beta} e_{\alpha}
\wedge e_{\beta},
$$
if $r_1^{\alpha\beta}$ and $r_1^{\gamma\delta}$ are
both not zero,
and $\alpha + \gamma $ is a root,
then $\beta + \delta$ must also be a root.
The last two terms in (2.19) are thus
$$
\sum
 r_1^{\alpha\beta}r_1^{\beta\delta}\bigl (e_{\alpha + \beta}
 \wedge e_{-\beta} \wedge e_{-\delta} +
e_{-\beta -\delta} \wedge e_{\alpha} \wedge e_{\gamma}
\bigr)
= d  \bigl(  \sum r_{22}^{\alpha\beta} e_{\alpha}
\wedge e_{\beta}\bigr)
$$
with
$$
r_{22}^{\alpha + \gamma , -\beta - \delta } =
- r_1^{\alpha\beta}r_1^{\gamma\delta}.\eqno(2.22)
$$

\q The calculation can be carried on to higher orders, but the
panorama of Belavin
and Drinfeld [1] is already emerging and it may be sufficient to
give the final answer.

\b
{\bf Proposition.} The exact deformations
 are classified up to equivalence by
$$r(\hbar) = r + \sum_{i = 1}^l \hbar^i \sum_{\sigma(i)}
  e_{\alpha}
\wedge e_{-\beta},
$$
where $\sigma = \sum\sigma(i)$ is the subset of $ \{(\alpha,\beta) \}$
 that satisfies the following additional conditions:  1) The first
entry ($\alpha$) runs over a subalgebra $\G_1$ of positive roots
  and the second
entry ($\beta$) runs over a subalgebra $\G_2$ of positive roots.
 2) Both subalgebras are generated by simple roots and
there is an isomorpism
$\tau: \G_1 \rightarrow \G_2$, such that some iteration $\tau^k$ of $\tau$
leads out of $\G_1; \,\, \tau^k \alpha \notin \G_1$ for all $\alpha \in \G_1$,
and $(\alpha, \beta) \in \sigma$ iff $\beta = \tau^m \alpha$ for some
positive integer $m$. 3) $w(\alpha) + w(\tau \alpha) = 0.$
 Finally, $(\alpha,\beta) \in \sigma(i)$ means
that both roots have the same length $i$, and $l$ is the rank of \g.

\r ?
\ve

{\ce{\bf 3. Cohomology.}}

\b

A. LIE ALGEBRAS.%3A

\q Let { \got g} be a Lie algebra over \Cit,
with structure tensor $\epsilon$,
and {\got g}* its vector space dual. The cochains $C_p^q$ are functions from
{\got g}$^{\wedge p}$ to {\got g}$^{\wedge q}$, and the action of $d$ on
$\sigma \in C_p^q$ is given by
$$
\eqalign{
 d\sigma (u_o,u_1,...,u_p) =
&\sum_i (-)^iu_i\sigma(u_0,...,\hat u_i,... ,u_p)\cr
&+ \sum_{i<j}(-)^{i+j}\sigma([u_i,u_j],u_o,..,\hat u_i,..,\hat u_j,..,u_p)\cr
}\eqno(3.1)
$$
In the first term the element $u_i \in \,\,${\got g} acts on
$\sigma$ according by the adjoint action:
$$
u\sigma := [\Delta_0^{p-1}u,\sigma] = [u,\sigma].\eqno(3.2)
$$
The second formula should be taken to be short hand for the first one, where
$$
\Delta_0u = 1 \otimes u + u \otimes 1,\eqno(3.3)
$$
and $\Delta^{p}$ is the $p$'th iterate of this primitive tensor product.
In view of what happens later, it is worth notice that this notion
of tensor product is needed from the beginning.

\q Let us examine the meaning of the simplest examples.
\b
\underbar {$p = 0.$} \enskip If $\sigma \in C_0^q$, then
$$
d\sigma(u) = [u,\sigma].\eqno(3.4)
$$
Thus $d\sigma = 0$ means that $\sigma$ is an
{\bf invariant element} of {\got g}$^{\w q}$, and the coboundaries
 $B_1^q \subset C_1^q$ are
those cochains that can be expressed as a commutator; in particular,
the {\bf derived algebra} is $ B_1^1($\g$) = [${\got g,g}$]$.   The function $C
= d\sigma$ from {\got g} to {\got g}
generates an infinitesimal {\bf inner homomorphism} of the Lie algebra,
$$
u \mapsto u_{\sigma} = u + \hbar [\sigma,u].\eqno(3.5)
$$
Consider the case $q = 2$;
$r \in C_0^2$ is an element of {\got g $\otimes\,\, $g}.
The solutions of $dr = 0$ are the invariant elements in {\got g $\otimes $ g},
and
$
f = dr
$
is a function from {\got g } to {\got g $\wedge $ g};
this special case is central to our subject.
\b
\underbar {$p = 1.$}\enskip
If $C\in C_1^q$, then
$$
dC(u,v) = [u,C(v)] - [v,C(a)] - C([u,v]).\eqno(3.6)
$$
If $q = 0$ then $C$ takes its values in   \Cit,
the first two terms disappear,and $dC = 0$ means
that $C$ is an invariant \Cit-valued linear functional on {\got g}.
If $q = 1$, then $C \in C_1^1$ is a function from {\got g} to {\got g};
this finds application as a generator of a change of basis in
the Lie algebra. Thus, consider the mapping
$$
u \mapsto u_c = u + \hbar C(u);\eqno(3.7)
$$
then to order $\hbar$
$$
[u_c,v_c] = [u,v] + \hbar ( [u,C(v)] - [v,C(u)]) =
   [u,v]_c + \hbar dC(u,v).\eqno(3.8)
$$
That is, $C$ generates a homomorphism of {\got g} iff $dC = 0$.
If $C = d\sigma$, then the homomorphism is an inner one, and
$dC = 0$ is automatic since $d^2 = 0$. The functions
$dC$ are the {\bf infinitesimal deformations} of the structure tensor induced
by a reparameterization of the Lie algebra.

\q The case $q = 2$: let $f$ be a function from {\got g} to {\got g $\w $ g},
$f$ can be interpreted as a {\bf coproduct},
$$
\Delta_f(u) = f(u).\eqno(3.9)
$$
If $df = 0$, then this coproduct is compatible with the
Lie structure of {\got g} and turns {\got g} into a {\bf Lie bialgebra}, and
if $f = dr$  we have a ``coboundary Lie bialgebra";
we postpone further discussion of this point.
\b

\underbar{$p = 2.$} \enskip  When $E \in C_2^0$ we
have $dE = \sum_{cyclic}E([u,v],w)$,
and $dE = 0$ is the condition for $E$ to define a {\bf central extension} of
\g.
 The Lie structure tensor $\epsilon$ belongs to $C_2^1$. For any
such cochain our formula (3.1) gives
$$
dE(u,v,w) = \sum_{cyclic}([u,E(v,w)] + E([u,v],w)\eqno(3.10)
$$
To understand the meaning of this condition consider a
deformation of the structure
$$
\epsilon \mapsto \epsilon' = \epsilon + \hbar E.\eqno(3.11)
$$
where  $\hbar $ is the deformation parameter.
It is easy to see that $dE = 0$ is
the condition that the deformed structure tensor satisfy the
Jacobi identity, to first order in $\hbar $. If $E = dC$, then the deformation
is
induced by a change of basis and $dE = 0$ trivially.
When the formula is applied to $\epsilon\,\,$  itself
it yields
$$
d\epsilon (u,v,w) = 2\sum_{cyclic}\epsilon([u,v],w) =
2\sum_{cyclic}\epsilon(\epsilon(u,v),w),\eqno(3.12)
$$
and $d\epsilon = 0$ is the Jacobi identity.
Thus (3.10) is the linearization of (3.12).

\b
\q Let $Z_p^q \subset C_p^q$ be the closed cocycles,
$B_p^q \subset Z_p^q$ the coboundaries, and $H_p^q = Z_p^q/B_p^q$.
We have seen that $H_0^q$ consists of the invariants in {\got g$^{\wedge q}$},
$H_1^0$ is the set of invariant functionals on {\got g},
$H_1^1$ is the set of outer automorphisms of {\got g}, $H_1^2$
classifies the non-coboundary Lie bialgebras and $H_2^1$
contains the essential deformations of the structure.
The obstructions to carrying the deformations to higher order
lie in $H_3^1$, and $H_2^0$
classifies the central extensions of \g.
The Haar measure (on the group) belongs to $H_n^0,\, n = $dim(\g).

\q If {\got g} is {\bf simple}, then $H_p^q = 0$ whenever
$q \neq 0$;  $H_1^0, H_2^0 $ also vanish, and $H_3^0 \neq 0$.
In particular, $H_1^2 = 0$
means that all structures of Lie bialgebra are of the coboundary type,
 $H_2^1 = 0$ means that the Lie structure is {\bf rigid}, and $H_2^0 = 0$ that
there are no central extensions.

\q Let $S \subset T$ be {\g}-modules and $Q$ the \g-module $T/S$; then
$T$ is called an {\bf extension} of $S$ by $Q$. As vector spaces,
$T = Q \oplus S$; let $\pi_{_0} = \pi_{_Q} \oplus \pi_{_S}$, then
there is a map $C$ from \g \enskip to Hom$(Q,S)$ such that
$$
\pi_{_T}(u) = \pi_{_0}(u) + C(u),
$$
where for $u \in $ \g,
$$
C(u): \,\,q \oplus s \mapsto 0 \oplus C(u)q.\eqno(3.13)
$$
The fact that $\pi_{_T}$ is a representation of \g\enskip amounts to $dC = 0$,
$$
[\pi_{_T}(u),\pi_{_T}(v)] - \pi_{_T}([u,v]) = [u,C(v)] - [v,C(u)] - C([u,v])
= dC(u,v).
$$
Suppose there is a map $\sigma : Q \rightarrow S$ such that
$$
C(u) = [\pi_{_0}(u),\sigma] = \pi_{_S}(u)\sigma - \sigma\pi_{_Q}(u)
= d\sigma(u);
\eqno(3.14)
$$
then $dC = 0$. We show that in this case the extension is trivial (split).
 Change basis by setting
$$
q \oplus s' = q \oplus (s + \sigma q),
$$
then using (3.14) one obtains
$$
\pi_{_T}: q \oplus s' \mapsto \pi_{_Q}q \oplus \pi_{_S}s'.
$$
Conclusion: if $C = d\sigma$, then $\pi_{_T}$ is equivalent to $\pi_{_0}$.
This example of cohomology is fundamental; some of the examples discussed
above can be interpreted this way, with proper identification of the spaces.

\bb

B. LIE BIALGEBRAS.%3B

\q Suppose now that \g$^*$ also has a Lie structure,
 with structure tensor $f$. By pullback, or duality, the cochains that
were seen as functions from \g$^{ \wedge p}$ to \g$^{\wedge q}$,
can be interpreted as functions from \g$^{*\wedge q}$ to \g$^{*\wedge p}$,
or more symmetrically,  from
\g$^{\wedge p} \,\otimes \,$\g$^{*\wedge q}$    to \Crm.
(In this last case it is natural to replace \enskip \Crm \enskip
by any space that is both a \g-module and a \g$^*$-module.)
Now we have two differential operators,
$$
d: C_p^q \rightarrow C_{p+1}^q, \q \partial: C_p^q \rightarrow C_p^{q+1}.
$$
This is what is called a {\bf bi-complex}, provided that $d, \partial$
anticommute:
$d \circ \partial + \partial \circ d = 0$. The structure is defined by
$\epsilon \in C_2^1$ and $f \in C_1^2$; the Jacobi identities are
$d\epsilon = 0,\,\, \partial f = 0$.
\b
{\bf Theorem.} The following statements are equivalent
$$
1) \,\,d\partial + \partial d = 0,\,\, 2) \,\,df = 0,
\,\,3)\,\, \partial \epsilon = 0.
$$
\b
Strictly speaking, an overall sign has to be adjusted in the definition
of $\partial$ to get the operators to anticommute.
The equivalence of 2) and 3) is easy to prove. On the dual bases
$(L_i,\Gamma^i)i = 1,...,$dim(\g),
$$
(df)_{ij}^{kl} = \epsilon_{im}^kf_j^{ml} + \epsilon_{im}^lf_j^{km} - (i,j) -
\epsilon_{ij}^mf_m^{kl} = (\partial \epsilon)_{ij}^{kl}.\eqno(3.17)
$$
The last equality (up to a sign, perhaps) is an evident consequence
of the symmetry of
the tensor with respect to the roles played in it by the two structure tensors.

\q When $df = 0$ we say that the structures $f$ and
$\epsilon$ are {\bf compatible}. The theorem gives
a first interpretation of what this means.
 From now on we are given a pair of compatible forms $f$ and $\epsilon$.

\q From the discussion of the connection between $d\epsilon = 0$
and the Jacobi identity it follows that the form $f$ is $\partial$-closed, so
$$
\partial f = 0, \q df = 0.
$$
The first equation says that $f$ defines a Lie structure on \g$^*$
and the second that it is compatible with the Lie structure of \g.

\q Alternatively, one can ascribe both structures to \g, associating
to $f$ the coproduct
 $\Delta_f$: \g $\, \rightarrow$ \g $\,\wedge \,$\g,
$$
\Delta_f(u) = f(u).
$$
The subscript is needed to distinguish $\Delta_f$ from the primitive
coproduct $\Delta_0$, Eq.(3.3). This double
structure makes \g \enskip into a Lie bialgebra:
$$
[u,v] = \epsilon(u,v), \q \Delta_f(u) = f(u), \q \partial f = df =
0 =\partial \epsilon = d\epsilon.
$$
In particular, if $f$ is a coboundary,
$$
f = dr, \q r \in C_0^2,
$$
then $df = 0$ automatically,
while $\partial f = 0$ becomes a condition on $r$,
$$
 \partial f = \partial dr = -d\partial r = 0.\eqno(3.18)
$$
That is, $\partial r$ must be an invariant element of \g$^{\wedge 3}$.
The formula is
$$
\partial r = [r_{12},r_{13}] + [r_{12},r_{23}] + [r_{13},r_{23}] =:
{\rm C\hskip-0.4mmY\hskip-0.7mmB(r)},
\eqno(3.19)
$$
and $\partial r = 0$ is the classical Yang-Baxter relation.

\q Strictly, this makes sense if $r$ is antisymmetric. It turns out that,
by including in $r$ a symmetric, invariant tensor, one can replace
$d\partial r = 0$ by the more convenient condition $\partial r = 0$.
{}From now we assume, following Drinfeld, that
$$
 r + r^t = K, \q \partial r = 0.\eqno(3.20)
$$
Here $K$ is the Killing form of \ge and $\partial r $ is defined by (3.19).
\ve
\q The dual, \g$^*$, also becomes a Lie bialgebra, with
$$
\{a,b\} = f(a,b), \q \Delta_{\epsilon}a = \epsilon(a).
$$
It is natural to ask if this is also of the boundary type.

\q Let $A,U$ be the matrices of the  adjoint actions of
$a \in \,$\g$^*$  \enskip and $u \in \,$\g,
$$
A_j^i = a_kf^{ik}_j, \q U_j^i = u^k\epsilon^i_{jk}.\eqno(3.21)
$$
The following are identities:
$$
\eqa{
&rU + Ur + f(u) = 0,\cr
&rA + Ar + r\epsilon(a)r = 0.\cr
}\eqno(3.22)
$$
The first is just a way to write $f(u) = [u,r]$, or $ f = dr$;
the second is the same as the classical Yang-Baxter relation, $\partial r = 0$.
If $r$, considered as a map from \g \enskip to \g$^*$, were invertible,
then it would follow that
$
\epsilon(a) = [a,r^{-1}],
$
or $\epsilon = \partial r^{-1}$. In fact, $r$ is not invertible,
nevertheless we shall see that

\r {NB}
%there is a map $r^*$ from \g$^*\,$ to \g,
%such that $\epsilon = \partial r^*$. This does not mean that \g$^*$ is a
%coboundary Lie bialgebra, since the map
%$r^*$: \g $ \rightarrow $\g  \enskip is not
%an element of \g$^* \otimes \,$\g$^*$.
%\b
%{\bf Theorem.} If \g \enskip is a coboundary Lie bialgebra,
%then the structure tensor  $\epsilon$ is the differential of a function
%$r^*$: \g $ \rightarrow $\g, $\epsilon = dr^*$.
%\b
%It means that the dual of \g \enskip is a sub-bialgebra
%of a coboundary Lie bialgebra. The proof is given below.
\bb
C. MANIN ALGEBRA.%3C

\q Another implication of compatibility of $f,\epsilon$
is the existence of a structure of Lie algebra on \g $\,\oplus \,$\g$^*$.
Given the brackets on \g \enskip and on \g$^*$ one adds
($a \in \,$\g$^*$   and $u \in \, $\g )
$$
[a,u] = {\rm coad}_au - {\rm coad}_ua,
$$
where
$$
{\rm coad}_au(b) = -u(\{a,b\}), \q {\rm coad}_ua(v) = -a([u,v]).
$$
The Jacobi identity is satisfied by virtue of $df = 0$. The Lie algebra
obtained in this way is denoted
\g $\, \triangleright \hskip-0.5mm \triangleleft \,$\g$^*$.
\b
{\bf Theorem.} The Manin algebra \g$\,\bow \,$\g$^*$ \enskip is
isomorphic to \g $\,\,\oplus $ \g.
\b

We give a constructive proof.

\q Consider the adjoint action of \g $\,\bow   $ \g$^*$ \enskip on itself.
Let    $u,v \in $ \ge and $a,b \in $ \g$^*$ \enskip; then
$$
ad_{(u,a)}(v,b) = \big ( vA - f(u)b - vU, Ub - Ab - \epsilon(a)v \big ).
\eqno(3.23)
$$
with $A,U$ as in (3.21).
The action of \g \enskip involves an extension of the coadjoint representation
by the adjoint representation. Since $f = dr$, we can trivialize this extension
by a change of variables. Define $v' = v - rb$, then
$$
{\rm ad}_{(u,0)}(v',b) = (-v'U,Ub),\eqno(3.24)
$$
by virtue of the first of Eq.s (3.22).
The action ad$_{(0,a)}$ likewise involves an extension,
but since \g$^*$   is not simple it is not clear whether
it is trivial.

\q The problem is to determine the invariant subspaces for the action
in (3.23). By (3.24) every subspace that is invariant
under ad$_{(u,0)}$ has the form
$$
S_{\kappa} = \{(v,b); v' = \kappa Kb\},
$$
with $\kappa  \in \,\,$\Cit \enskip and $K$ the Killing form.
This subspace is invariant under ad$_{(0,a)}$ iff
$$
A(\kappa K-r) + (\kappa K-r)A - (\kappa K-r)\epsilon(a)(\kappa K-r) = 0.
$$
When $\kappa = 0$ this reduces to the second of Eq.s(3.22), which is
the classical Yang-Baxter relation;
the latter is also obeyed by $r^t = K - r$,
so the statement is valid for $\kappa  = 1$ as well.
That is; $S_0,S_1$ are invariant subspaces of the adjoint action of
\g $\, \bow \,$\g$^*$. The actions are determined by
$$
{\rm ad}_{(u,a)}: b \mapsto
 \matrix{(U - A - \epsilon (a)r)b,&(v,b) \in S_0, \cr
  (U - A + \epsilon(a)r^t)b, & (v,b) \in S_1.\cr}
$$
The two subalgebras commute, $[S_0,S_1] = 0$; $(ra,a) \in S_0$ and
$(r^ta,-a) \in S_1$ act by the matrix $-$ad$_{Ka}$ on $S_0$ and on
$S_1$, respctively. The theorem is proved,
 \g $\,\bow \,$\g$^*$ \enskip is isomorphic to
$S_0 \oplus S_1 = \,\,$\g $\,\oplus \,$\g. This result was first obtained
by Reshetikhin and Semenov-Tian-Shanskii [14].
\ve

\q Consider now the bicomplex

$$
\diagram{
C^0_0 & \harr{d}{} & C^0_1 &
\harr{d}{} & C^0_2 & \harr{d}{} & C^0_3
& \ldots\cr
\varr{\partial}{}&&\varr{\partial}{}&&\varr{\partial}{}&&
\varr{\partial}{}&\cr
C^1_0 &\harr{d}{}&C^1_1 &
\harr{d}{}&C^1_2 &\harr{d}{}&C^1_3 &
\ldots\cr
\varr{\partial}{}&&\varr{\partial}{}&&\varr{\partial}{}&&
\varr{\partial}{}&\cr
C^2_0 &\harr{d}{}&C^2_1 &
\harr{d}{}&C^2_2 &\harr{d}{}&C^2_3 &
\ldots\cr
\vdots && \vdots && \vdots &&\vdots &\cr
}
$$
{}From the perspective of \gg$^*$, \g \enskip and \g$^*$   are
subalgebras, and the chains are
$$
C^{i-2} = C^i_0   \,\,\cup\, C^{i-1}_1 \,\,\cup ...
\cup \, C^1_{i-1} \,\,\cup \, C^0_i,
$$
functions from (\gg$^*$)$^{\w n}$ to  \thinspace  \Crm \enskip
(or to some \gg$^*$-module M).
Thus it is
natural to consider the chain
$C^1 = C^2_1 \,\,\cup\,C^1_2$, where we find the
tensors $f$ and $\epsilon$ united,
both being interpreted as functions from \gg$^* \,\,\otimes$ \gg$^*$ \enskip to
 \gg$^*$.
It is customary to start the complex at $p,q = 1,1$;
that explains the strange notation. The cohomology space associated with
$C^1$ is called $H^1$ or more precisely $H^1$(( \gg$^*$ )$^{\w 2}$) since
$H^1(M)$
normally is taken to mean the homology group with values in the module M.

\q The complex built on these chains is called the total complex.
The differential operator of the total complex is $d + \partial$.
It is natural to use this complex to investigate the deformations
of \gg$^*$; that is, simultaneous deformations of $f$ and $\epsilon$.
Since \g \enskip is simple, \gg$^*$ \enskip is semisimple and $H^1(M) = 0$
for all $M$. In particular, that $H^1$((\gg$^*)^{\w 2}$) is trivial
means that \gg$^*$ \enskip is rigid, all Lie
algebra deformations are trivial.
Of course, this does not mean that the deformations are uninteresting.
We saw that the classical Yang-Baxter relation
means that $\partial r = 0$, while $dr = f $
does not vanish; the classification of r-matrices is
not a problem of cohomology on the Manin algebra.

\ve

{\ce{\bf 4. Quantization.}}
\b
A. SYMPLECTIC STRUCTURE AND QUANTIZATION.%4A

\q
 Let \g \enskip be a simple, complex Lie algebra,
\g$^*$ its vector space dual, $G$ the associated local
Lie group and $A_0 = A_0(G)$ the germ of functions
at the identity. Whenever convenient we refer to  dual bases
$(L_i, \Gamma^i)i = 1,2,...,n$ for \g  \enskip and \g$^*$.
The $\Gamma^i$ may be interpreted as a set of exponential coordinates on $G$;
for $g \in G$ one has $g = e^{\Gamma(g)}, \,\Gamma (g) =
\Gamma^i(g)L_i$. The structure tensor for \ge is denoted $\epsilon$,
$[L_i,L_j] = \epsilon_{ij}^kL_k$.

\b

\q  The view of quantization presented here
is the same as that of ref [11],
with the new elements introduced by Drinfeld [2]. One regards $G$ as a
phase space  and $A_0(G)$ as the space of classical observables.
We suppose that $A_0$ is endowed with a Poisson bracket $\{ , \}$,
determined by a Poisson form $f$,
$$
\{\Gamma^i,\Gamma^j\} = f_k^{ij}\,\Gamma^k,\eqno(4.1)
$$

The coefficients are antisymmetric in $i,j$ and satisfy the Jacobi identity,
$$
\sum_{(ijk)}f_n^{im}\,f_m^{jk} = 0.\eqno(4.2)
$$
It is possible to identify phase space with a symplectic leaf of $G$,
rather than $G$ itself.
\q The space $A_0(G)$ is thus furnished with two independent structures; the
commutative algebra of functions,
$$
(ab)(g) = a(g)b(g),\q a,b \in A_0,\q g \in G,\eqno(4.3)
$$
and the Poisson-Lie bracket [2]
$$
\{a,b\} = \Lambda^{ij}\bigl[(L_ia)(L_jb) - (L_i'a),(L_j'b)\bigr], \q
\Lambda^{ij} = \Gamma^k\,f_k^{ij},\eqno(4.4)
$$
where $L_i(L_i')$ is to be understood as a
vector field of left (right) translations.

\q The aim of quantization is to replace both structures by a single,
noncommutative product denoted $a*b,\,\, a,b \in A_0(G)$,
expressed as a formal deformation of the Poisson bracket,
$$
a*b = \sum_{n = 0}^{\infty}\hbar^nC_n(a,b) =
ab + (i\hbar/2)\{a,b\} + o(\hbar^2),\eqno(4.5)
$$
with deformation parameter $\hbar$ (Planck's constant).

\q We insist on associativity of the new product,
and this implies conditions on the $C_n{}'s$ that are best formulated
in terms of the Hochschild cohomology of $A_0$.

\q The Hochschild cochains  are in this instance linear functions
from $A^{\otimes n} $ to $A$ and
$$
\eqalign{
dC(a_1,...,a_{p+1}) = &a_1  C(a_2,...,a_{p+1})\cr
&+ \sum(-)^iC(a_1,..a_ia_{i+1},...a_{p+1})\cr
&\q + (-)^{p+1}C(a_1,...,a_p)a_{p+1}.
}
$$

\q The associativity constraint is
$$
\eqalign{&
a * (b*c) - (a*b)*c = \sum_0^{\infty} \hbar^nD_n(a,b,c) = 0,\cr
&
D_n = E_n - dC_n \in C_3^1,\cr
& E_n(a,b,c) = \sum_{r =1}^{n-1}
C_r(C_{n-r}(a,b),c) - C_r(a,C_{n-r}(b,c)).
\cr}
$$

One can verify that, if $D_n = 0$ for $n \leq m$ then $dE_{m+1} = 0$.
To continue the deformation to the next order one must find $C_{m+1}$
such that $E_{m+1} = dC_{m+1}$, Therefore $E_{m+1}$ must be in the zero class
of $H_3^1$ ( usual name $H^3(A)$). One can also show that two $*$-products are
equivalent if
the difference of their cochains are exact at each level; inequivalent
$*$-products are therefore classified by $H_2^1$ ($ = H^2(A)$).

\q The odd part $a*b - b*a$ of the product defines a structure of Lie
algebra. The Jacobi identity implies conditions on the $C_n{}'s$
that are expressed in terms of the Chevalley-Eilenberg
cohomology of $\{,\}$
and   are automatically satisfied by associativity of $*$.

 \q We denote by $A_{\hbar}(G)$ the algebra of
  formal series in $\hbar$ with coefficients in $A_0$,
 with an associative $*$-product.
\bb
B. $*$-PRODUCTS AND ABSTRACT ALGEBRAS.%4B

\q
Quantization as exposed so far is characterized by a well
defined correspondence  between the quantum algebra $A_{\hbar}$
and the classical algebra $A_0$; every element in $A_{\hbar}$ is identified
with a  power series in $\hbar$ with coefficients in $A_0$
and, when the series converges, with a unique function on phase space.
 The advantage of this approach is considerable [11],[12];
for example,  integration on the quantum algebra is readily available.

\q A more abstract point of view is to study the abstract, unital
algebra generated by elements $(\Gamma_i)i = 1,2,...,n$,
with relations such as
$$
\Gamma^i\,*\,\Gamma^j - \Gamma^j\,*\,\Gamma^i =
i\hbar\,f_k^{ij}\,\Gamma^k. \eqno(4.6)
$$
Apart from the fact that this may be a proper subalgebra of
 $A_{\hbar}$ (which is not the issue I want to  discuss) I point out that one
here abandons
the actual interpretation of an element of the algebra $A$
as a
function on $G$, making the term ``quantum algebra of functions"  somewhat
misleading. (Relations of the type (4.6) do not necessarily hold,
the point is that when they do, then the alternative,
abstract point of view becomes available. This may also be possible with
other types of relations, as with certain types of quadratic relations
(Sklyanin algebras), but a detailed investigation is required to
determine the nature of an abstract algebra defined by relations   not
of the Lie type.)

\b

C. COMPATIBILITY.%4C

\q We have made no use, so far, of the group structure on
our phase space. It gives rise to a condition of ``compatibility"
between the Poisson structure on \g$^*$ and the Lie bracket on \g,
namely $df = 0$, or in components
$$
[L_i,f_j] - [L_j,f_i] - \epsilon_{ij}^kf_k = 0, \q
f_k := f_k^{mn}L_m \otimes L_n.\eqno(4.7)
$$
The real meaning of this relation must be clearly grasped; we begin by
eliminating some potential misconceptions. See Drinfeld, ref. [2].

\q The Lie structure   of \g$^*$ gives rise to a coproduct on
\g,
$
\Delta_f:$ \g $\,\rightarrow$ \g $\,\,\otimes$ \g ,
 defined by
$$
\q  \Delta_f(L_k) = f_k.  \eqno(4.8)
$$
This is a standard application of duality,
$$
(\Delta_fL_k)(\Gamma^m,\Gamma^n) = L_k(\{\Gamma^m,\Gamma^n\}).
$$
Satisfying the condition of ``compatibility" does not turn $\Delta_f$
into a homomorphism, neither does $\D_f$ become
an intertwiner the adjoint action, though
this misses by a factor of two.

\q A closer look at  (4.7) reveals that it involves the
standard coproduct on the
enveloping algebre $U$(\g), generated by
$$
\Delta_0: L_i \mapsto L_i \otimes 1 + 1 \otimes L_i,\eqno(4.9)
$$
since the first two terms in (4.7) involves the adjoint action of
 \g \enskip on \g $\,\otimes\,$\g.
So the real meaning of the   compatibility condition
should be sought in the context of the
enveloping algebra.

\q Indeed,  the condition of compatibility, Eq(4.7),
is the condition that there is a homomorphism to order $\hbar$,
($U = U$(\g) from now on)
$$
\Delta_{\hbar}: U  \rightarrow U  \otimes U
$$
generated by
$$
\Delta_{\hbar}(L_i) = L_i \otimes 1 + 1 \otimes L_i +
(i\hbar/2)\,f_i^{mn}\,L_m \otimes L_n; \eqno(4.10)
$$
that is,
$$
\Delta_{\hbar}([L_i,L_j]) = [\Delta_{\hbar}(L_i),\Delta_{\hbar}(L_j)] +
o(\hbar^2).
$$
This last coincides with (4.7) and expresses the true meaning of
the compatibility condition.
This interpretation is dual to Drinfeld's Poisson Lie groups [2].

\q This establishes a firm link between the compatibility condition,
 with its well
known connection to bialgebra structure and cohomology,
 and the type of deformation
that characterizes quantization.
The condition is the new element that comes into
play when phase space is a group manifold.

\bb
D. FURTHER COMMENTS ON QUANTIZATION.%4D

\q
Let a $*$-product be given on $A_0(G)$ and suppose for concreteness
that relations (4.7),
$$
\Gamma^i\,*\,\Gamma^j - \Gamma^j\,*\,\Gamma^i =
i\hbar\,f_k^{ij}\,\Gamma^k. \eqno(4.11)
$$
hold exactly. As we said, we may choose to forget the
origin of this relation and   use it only to define an abstract algebra
generated by elements $(\Gamma^i)i = 1,2,...,n$. Now conversely,
given this abstract algebra, we may try to realize it concretely
in terms of a $*$-product on $A(G)$. This amounts to establishing
(choosing) a mapping from the abstract algebra into $A_0$;
such a map is called an ``ordering".
The most traditional way to do this is to set up a basis
for the abstract algebra, for example
$$
(\Gamma^1*)^{k_1}...(\Gamma^n*)^{k_n} =: \Gamma_*^{\{k\}},
$$
and a linear correspondance that associates each such to
a specific function on $G$, for example, to the function
$$
(\Gamma^1)^{k_1}...(\Gamma^n)^{k_n} =: \Gamma^{\{k\}},\eqno(4.12)
$$
Normal ordering and standard ordering are both of this type.
The same method can be applied when the relations are of quadratic type.

\q Our favorite method is a generalization of symmetric ordering. Let
$$
e_*^{\Gamma\cdot L/i\hbar} := \sum_{n = 1}^{\infty}{1 \over n!}
\left({1\over i\hbar}\Gamma\cdot L*\right)^n, \q
\Gamma \cdot L := \Gamma^iL_i ;\eqno(4.13)
$$
where the $L_i$ are interpreted as coordinates
on \g$^*$   and on $G^*$, the formal group generated by \g$^*$.
To this expression,   interpreted
as a formal series, one associates a specific function $E_*(L_1,...,L_n)$ on
$G$,
for example the function
$$
e^{\Gamma\cdot L/i\hbar} := \sum_{n = 1}^{\infty}{1\over n!}
\left({1\over i\hbar}\Gamma\cdot L\right)^n.
$$
It should be emphasized that the choice that is made for this function $E_*$
 (the ``star-exponential") is far from being without consequence; the image
of the ordering map may even be finite dimensional for certain choices.
 See ref. [15].

\q The expression (4.13) can also be interpreted as an element of the
local group $G^*$ with
Lie algebra \g$^*$ and structure tensor $f$.
We have
$$
 e_*^{\Gamma\cdot L/i\hbar} \otimes e_*^{\Gamma\cdot L/i\hbar} =
e_*^{\Gamma\cdot C(L)/i\hbar},\eqno(4.14)
$$
where
$$
C(L_i) =  L_i \otimes 1 + 1 \otimes L_i +
(i\hbar/2)\,f_i^{mn}\,L_m \otimes L_n + 0(\hbar^2); \eqno(4.15)
$$
This is just the Campbell-Hausdorff formula; in it $L_i\otimes L_j$
stands for coordinates for \g$^* \otimes \,$\g$^*$. To order $\hbar$,
$C(L)$ coincides with $\Delta_{\hbar}(L)$, so the following conjecture
$$
 e_*^{\Gamma\cdot L/i\hbar} \otimes e_*^{\Gamma\cdot L/i\hbar} =
e_*^{\Gamma\cdot \Delta_{\hbar}(L)/i\hbar},\eqno(4.16)
$$
yields $\Delta_{\hbar}$ correctly to order $\hbar$. In fact,
if \g \enskip were abelian, then this formula would
define a coproduct on $U $, compatible with the (trivial) Lie structure
of \g.

\q We can reverse the roles of \g \enskip and \g$^*$ (of $\epsilon$ and $f$),
since the compatibility condition is symmetric, and consider
$$
 e_*^{\Gamma\cdot L/i\hbar} \otimes e_*^{\Gamma\cdot L/i\hbar} =
e_*^{ \Delta_{\hbar}(\Gamma)\cdot L/i\hbar},\eqno(4.17)
$$
with
$$
\Delta_{\hbar}(\Gamma^k) = \Gamma^k \otimes 1 + 1 \otimes \Gamma^k
+ (1/i\hbar)\,\epsilon_{ij}^k\,\Gamma^i \otimes \Gamma^j.\eqno(4.18)
$$
This brings to mind the universal T-matrix [9],[10].

\bb
E. THE OPERATORS $d$ AND $\partial$ ON $U_{\hbar}$.%4E

\q
 Denote by $U_{\hbar}$ the bialgebra dual of $A_{\hbar}(\g^*)$,
endowed with a (possibly deformed) algebraic structure and the
coproduct $\Delta_{\hbar}$ as in Section 4D.
Compatibility is generalized or lifted up in the most natural way
from (4.9) to
$$
\Delta_{\hbar}(uv) = \Delta_{\hbar}(u)\,\Delta_{\hbar}(v),
\q u,v \in U_{\g}.
$$
We continue to use a $*$ for the product on $A_{\hbar}$, whenever that serves
to facilitate the interpretation of the formulas, but not for the
product on $U_{\hbar}$,
and write $U,A$ and $\D$ without the subscript from now on.
Unless explicit exception is made, all that follows is independent of the
(possibly restrictive) condition that $\D$ be expressible in the
form (4.18).

\q The basic cohomology on $A$ is Hochschild, with cochains
$$
C^q_p \in {\rm Hom}(A^{\otimes p}, A^{\otimes q}),
$$

   The formula for the differential operator is
$$
\eqalign{
dC(a_1,...,a_{p+1}) = &\Delta^{q-1}(a_1)C(a_2,...,a_{p+1})\cr
&+ \sum(-)^iC(a_1,..a_ia_{i+1},...a_{p+1})\cr
&\q + (-)^{p+1}C(a_1,...,a_p)\Delta^{q-1}(a_{p+1}).
}\eqno(4.19)
$$
with $a_1,...,a_{p+1} \in A$, $C(a_1,...,a_p) \in A^{ \otimes q}$.
To get the flavor of it,
consider the case $p = 0, q = 1$, so that $C \in A$, then from (4.19),
the first and the last terms contributing,
$$
(dC)(a) = aC - Ca;
$$
and (4.19) can be considered as a generalization of this formula, in two steps.
The generalization to the case where $C$ is in $A^{\otimes q}$ needs an action
of
$a \in A$ on this object; it is given by $\D^{q-1}(a)$. The generalization
to higher forms follows the pattern set by the de Rham complex.

 \q  The dual $U$ of $A$ is an algebra, and it has its own Hochschild
complex with differential operator $\partial$. By duality, this operator
can be expressed in terms of the coproduct of   $A$.
 While $d$ increases the number of arguments,
leaving the value in the same space, $\partial$ maps   the function
taking values in $A^{\otimes q}$ to one taking values in $A^{\otimes q+1}$.
To write it we need some notation. The coproduct can be expressed as a sum,
$$
\Delta(a) = \sum_{(a)} a^1 \otimes a^2,
$$
and this will be the meaning of $a^1,a^2$ in what follows. The formula is
$$
\eqalign{
\partial C(a_1,...,a_p) = &\sum_{(a_1,..,a_p)} a_1^1a_2^1...a_p^1 \otimes
C(a_1^2,...,a_p^2)\cr
& \q + \sum_i (-)^i\Delta_i C(a_1,...a_p)\cr
&\q + (-)^{q+1}\sum_{(a_1,..,a_p)} C(a_1^1,...,a_p^1) \otimes a_1^2...a_p^2.
\cr}\eqno(4.20)
$$
This can be considered as a generalization of $\D: U \rightarrow U \otimes U$,
when $p = 0, q = 1, C \in U$ and $dC = \D C$.
To illustrate the formula, consider the first term in the case that $p = 1$,
then for $u \in U$,
 $$
(uC)(a) =   (u \otimes C)(\Delta a) = \sum_{(a)} u(a^1) C(a^2).\eqno(4.21)
$$
 Let us indicate the values $C(a;u_1,...,u_q)$ of $C(a)$.
The last formula becomes
$$
\eqalign{
(uC)(a;u_1,...,u_q) &= (u \otimes C(u_1,...,u_q))(\Delta a)
 = \sum_{(a)} u(a^1) C(a^2;u_1,...,u_q)\cr & = \bigl(
\sum_{(a)}a^1 \otimes C(a^2)\bigr)(u \otimes u_1 \otimes ... \otimes u_q).\cr
}
$$
We see that it is the value of a one-form on $A $;
it is the value of the first contribution to $\partial C(a)$ at arguments
$u,u_1,...,u_q$.

\q Let us go a step further, to the case $p = 2$,
by representing $C(a_1,a_2)$ as $(\Delta D)(a_1,a_2) = D(a_1*a_2)$.
(This is sufficient, by linearity.) Then one gets an action of $u \in U$,
from Hom$(A^{\otimes 2},A^{\otimes q})$
to the same with $q$ replaced by $q+1$,
$$
(uC)(a,b;u_1,...,u_q) = \sum_{(a)(b)}\D (u)(a^1,b^1) \otimes
C(a^2,b^2;u_1,...,u_q).
$$
In general
$$
(uC)(a_1,...a_p)(u_1,...,u_q) =
\bigl (\Delta^{p-1}(u) \otimes
 C(u_1,...u_q)(\Delta a_1,...,\Delta a_p) \bigr),\eqno(4.22 )
$$
with the understanding that $\Delta^{p-1}(u)$ takes the arguments
$a_1^1,...,a_p^1$ and $C(u_1,...,u_q)$
 the arguments $a_1^2,...,a_p^2$.
This is the value of the first term in Eq(4.20) at $u,u_1,...,u_q$.

\q Checking the properties $d^2 = 0$ and $\pa^2 = 0$,
we find that they are equivalent to
$$ \D(a)\D(b) = \D(a*b), \,\,{\rm and} \,\,(1 \otimes \D)\D = (\D \otimes1)\D,
$$
respectively.

\b

\q  As a first application, let us take $p = q = 1$,
and write the definitions again,
$$
\eqalign{
&dC(a_1,a_2;b) =  (a_1 C(a_2))(b) - C(a_1a_2)(b) +( C(a_1) a_2)(b),\cr
&\partial C(a;b_1,b_2) = \sum_{(a)} (a^1 \otimes C(a^2))(b_1,b_2)
- (\Delta C(a))(b_1,b_2)
+ \sum_{(a)} (C(a^1) \otimes a^2)(b_1,b_2),\cr}
$$
The most interesting cochains are  $\UT \in C_1^1$, defined by
$\UT(a) = a$, the multiplication map $m \in C_2^1$, and the
coproduct $\D \in C_1^2$:
$$
d\UT = m, \q \partial \UT  = \D .\eqno(4.23)
$$
It follows that $dm = 0 = \pa\D$, and indeed
$$
\eqalign{
dm(a,b,c) &= 2\bigl(  a*(b*c - (a*b)*c)\bigr),\cr
\pa \D& = 2\bigl((1 \otimes \D)\D - (\D \otimes 1)\D\bigr).}\eqno(4.24)
$$
In the case of a Hopf algebra the properties of the counit and
the antipode can also be expressed
neatly in terms of the differential operators.
\b
{\underbar {Exercise.}}
When the operator $d$ is applied to $C_p^0 = U^{\otimes p}$, it
yields a map from $U^{\otimes p}$ to $U^{\otimes p+1}$. So in
this case $d$ is represented by an element $C \in C_p^{p+1}(U)$.
Investigate the question of whether these cochains are
closed and/or coclosed.
\ve
F. THE R-MATRIX.%4F

\q In this section let $\D$ denote the coproduct on $U$ and
$\D' $ the opposite coproduct,
$$
\D(u) = \sum_{(u)} u^1 \otimes u^2 , \q \D'(u) = \sum_{(u)} u^2 \otimes u^1.
$$
 Suppose further that there is an element $R \in U \otimes U$ such that
$$
\D'(u)R = R\D(u).\eqno(4.25)
$$
We have $\D' = \sigma \D$, where $\sigma$ interchanges the factors,
so if we define $P \in U \otimes U$ by
$$
P = \sigma R,
$$
then Eq(4.25) becomes
$$
\D(u)P - P\D(u) = 0, \q {\rm or }\q \pa P = 0.\eqno(4.26)
$$
Remember that $\pa$ is the ``$d$" of the algebra $U$,
applied to our complex via the identification $u(a) = a(u)$.
Thus to apply the formula (4.20) for $\pa$ one should interpret $P$
as a function from $A \otimes A$ to \Crm.
Then $\pa P$ is a function from $A \otimes A$ to $A$, which attests to its
having to do with multiplication on $A$. By these conventions we
avoid having to deal with two different double complexes, each with two
differentiations.

\q The matrix $R$ is said to be ``unitary" if $R_{12}R_{21} = 1$.
We have $\sigma R_{12}R_{21} = P^2$, and evidently $dP^2 = 0$. This points
at a hierarchy of operators,
$$
I_1 = R, \,\, I_2 =
R_{12}R_{21}, \,\,
I_3 =  R_{12}R_{13}R_{23}/R_{23}R_{13}R_{12}, ...\, .
$$
In terms of them we can express triviality, $I_1 = 1$,
unitarity, $I_2 = 1$ and the
Yang-Baxter relation, $I_3 = 1$.
Multiplying by $\sigma$ one turns these operators into the invariants
$$
J_1 = P, \,\,J_2 = P^2,\,\,
J_3 = (P_{12}P_{23}P_{12}/P_{23}P_{12}P_{23}).
$$
We have seen that $\pa P = 0$ and thus $ \pa P^2 = 0$ as a direct consequence
of (4.24).
As to the third invariant there is a well known algorithm
using associativity of the product on $A$ that leads to $ \pa J_3 = 0$.
This last is weaker than the braid relation, which makes $J = 1$ and which
is equivalent to the Yang-Baxter relation. I am unaware of
any convincing argument that justifies YB on the basis of associativity alone.
The statement that $\pa J = 0$ is trivial since $\pa P = 0$; the fact that
it can be derived using associativity is  fortuitous.
\b
\q Before going on I want to rephrase the theory in the Woronowicz picture.
The local group $G$ (with Lie algebra \g)
is contained in the (undeformed) enveloping
algebra $U_0$ of \g. For $g \in G$ the standard
(classical) coproduct takes the form
$$
\D_0 (g) = g \otimes g.\eqno(4.27)
$$ Assume that there is
a subset $G$ of the deformed algebra $U$ that preserves the
structure of the group, and that the group algebra is dense in $U$.
(To simplify, we can even assume, with no
essential loss of generality [13], that the algebraic structure of $U$
is preserved; this is what is called a preferred deformation [16].)
Let $\pi$ be a fundamental representation of $G$, of dimension $n$,
and extend it to $U$.

\q Following Woronowicz [6], define the functions $(T_i^j)i,j = 1,...,n$
on $U$ by
$$
T_i^j(g) = (\pi g)_i^j, \q g \in G.\eqno(4.28)
$$
The coproduct on $A$ is generated by
$$
\eqalign{
 (\D T)_i^j(g_1 \otimes g_2) & = (\pi g_1g_2)_i^j = \pi( g_1)_i^k \,\pi
(g_2)_k^j,\cr
 \D T_i^j& = T_i^k \otimes T_k^j.\cr
}$$
Another application of duality gives
$$
(T_i^j * T_k^l)(g) = \D(\pi g)_{ik}^{jl},\eqno(4.29)
$$
where if $\D(g) = \sum_{(u)} u^1 \otimes u^2$, then
$\D(\pi g) = \sum_{(u)} \pi(u^1) \otimes\pi(u^2)$. Now evaluate Eq(4.26)
in $\pi$, replacing $u$ by $g$:
$$
\D(\pi g)P - P\D(\pi g) = 0,
$$
where $P$ now should be understood as the matrix $\pi P$.
In view of (4.29) this can be read as
$$
T_i^j * T_k^l\, P_{jl}^{mn} - P_{ik}^{jl}\,T_j^m * T_l^n = 0,
\q {\rm or }\,\,\,\, [T*T,P] = 0.\eqno(4.30)
$$
We have thus learned that these relations,
that define the structure of the algebra generated by
the $T_i^j$, is one way to write $\pa P = 0$.

\q The relations (4.29) have the form $Q = 0$, where $Q$ is
a fixed mapping from $U$ to $U \otimes U$ , or from $A \otimes A$ to $A$.
Since $Q = \pa P$ it is evident that $\pa Q = 0$.
More generally, suppose that the relations of $A$ have the form $Q = 0$,
with $Q$ an element of $C_2^1$. What do we buy with
the restriction $\pa Q = 0$?
Interpreting $Q$ as an element of $C_2^1$ we have
$$
(\partial Q)(a,b) = \sum_{(a)(b)} a^1*b^1 \otimes Q(a^2,b^2)
- \D Q(a,b) + \sum_{(a)(b)} Q(a^1,b^1)\otimes a^2*b^2.
$$
Evaluate this two-form at $(u,v) \in U \otimes U$.
Assuming that $Q(a,b)$ vanishes at $u$ and at $v$,
we obtain
$$
\D Q(a,b)(u,v) = Q(a,b)(uv) = 0;
$$
that is, the relations are compatible with the structure of $U$. In
other words, if $T_1,T_2$ are two sets of generators, of commuting
copies of $A$, then the matrix elements of the matrix product $T_1T_2$
satisfies the relations of $A$. This is just what makes $A$ into
a bialgebra; in the application to physical  models, where
$T$ is the transition matrix, it is the key to integrability.

\q Conclusion: If the relations of a bialgebra $A$ are given
in the form $Q = 0$, $Q$ an element of $C_2^1$, then $\pa Q = 0$;
among such bialgebras those characterized by an R-matrix are precisely those
for which $Q$ is exact. This justifies the terminology of Drinfeld,
who calls them coboundary bialgebras. To be precise, one should insist
that \underbar {all} the relations that define $A$ are
contained in the statement $Q = 0$; the restriction of quantum
$gl(n)$ to quantum $sl(n)$ by fixing the quantum determinant
would seem to place the latter outside the category of coboundary bialgebras.

\ve

\ce {\bf 5. Deformations of twisted quantum groups.}
\b
\q We are interested in the deformations of quantum groups;
that is, the deformations of certain coboundary bialgebras
 inside the category of coboundary bialgebras. We pose
$$
R(\ep) = R + \ep R_1\eqno(5.1)
$$
and require that this matrix satisfy the Yang-Baxter relation to order $\ep$.
This is a linear condition on $R_1$, but it is quadratic in $R$, which
makes it complicated; besides, the cohomological
interpretation is not clear to me. What follows is a search for an
alternative formulation, invoking some ideas
from $*$-quantization already utilized by Drinfeld.

\q Suppose there is an operator $F$ on $U \otimes U$ such that the (deformed)
coproduct $\D$ of $U$ can be expressed as
$$
\D = F  \, \D_0,
$$
where $\D_0$ is the standard coproduct (4.27). This is just the dual
image of the formula (4.5) of a $*$-product on $A$, where $F$ is
the bidifferential operator defined by
$
a*b = F(a,b).
$
The simplest possibility is that $F$ can be identified with an element of
$U \otimes U$. Attempts to construct such an element in the case of the
simplest
quantum group were frustrated, and indeed it does not exist, as shown
by K.M. Lau. His argument is this.  For $u,v \in U$,
$$\
F \, \D_0(uv) = \D (uv) = \D (u) \D (v) = F \, \D_0(u) \,F \,\D_0(v)
$$
If $F$ is invertible we get
$
\D_0(u)\D_0(v)  = \D_0(u)F\, \D_0(v)
$
and thus $F = 1$. Therefore, more generally let $F$ be an element
of $C_2^2$. According to ref. [13], there is no loss of generality in
taking $F = {\rm ad}\, N, \,N \in U \otimes U$, invertible,
$$
 \D (u) = {\rm ad}\, N \, \D_0(u) = N \,\D_0(u) \,N^{-1}.\eqno(5.2)
$$
Calculations on the simplest quantum group confirms that such $N$ exists.
Note also that this formula is suggested by (4.4).

\q Consider the P-matrix for $U_{qq'}(gl_2)$,
$$
P = \sum_iM_i^i \otimes M_i^i  + (1 - q'/q)M_1^1 \otimes M_2^2 +
q'M_1^2 \otimes M_2^1 + q^{-1}M_2^1\otimes M_1^2.\eqno(5.3)
 $$
(The P-matrix for twisted quantum $gl(n)$ is very similar and the calculation
of $N$ that follows generalizes immediately to that case.)
with parameters $q$ and $q'$. The associated starproduct is not defined, but
the relations among the matrix elements of
$$
T = \pmatrix{a&b\cr c&d \cr}
$$
 implied by (4.30) are
$$
\eqalign{
a*b &= qb*a,\,\, c*d = qd*c,\,\, c*a = q'a*c,\,\, d*b = q'b*d,\cr
c*b &= qq'b*c,\,\, a*d-d*a = (q-q')b*c.\cr}\eqno(5.4)
$$

We ask for a pair of matrices
$ N,\tilde N$ such that the product defined by
$$
TT =   \tilde N T*T N
$$
is abelian. In other words, the $*$-product would be
obtainable from an abelian
product by a twist
$$
T*T =  \tilde N^{-1} TT N^{-1},
$$
this being equivalent to (5.2) if $ \tilde  N = N^{-1}$.
The calculations can be done very efficiently with the aid
of quantum planes and antiplanes, covariant and contravariant;
the general solution is
$$
N = 1 \oplus \pmatrix{\alpha & \beta\cr \gamma & \delta \cr}
\oplus 1, \,\,
 \tilde  N = 1 \oplus \pmatrix{\tilde \alpha & \tilde \beta\cr
\tilde \gamma & \tilde \delta \cr}
\oplus 1,
$$
with
$$
\eqalign{
\pmatrix{\alpha & \beta\cr \gamma & \delta \cr} &=\,
\pmatrix{q+q' & q-q' \cr 0 & 2qq' \cr} \pmatrix{A &B\cr B&A}\cr
\pmatrix{\tilde \alpha & \tilde \beta\cr \tilde \gamma & \tilde \delta \cr
}& =\,
\pmatrix{\tilde A&\tilde B \cr \tilde B&\tilde A\cr}
\pmatrix{2qq'&q'-q \cr 0 & q+q'\cr}.\cr }\eqno(5.5)
$$
The coefficients $A,B$ are arbitrary except $A^2 - B^2 \neq 0$.
The basis  used in $V \otimes V$ is $(11,12,21,22)$,
 so for example $\beta = N_{12}^{21}$.
The ``opposite product" (dual to $\D'$) is characterized by similar relations,
in which $q,q'$ are replaced by $1/q,1/q'$.
The general expressions for the associated N-matrices are thus
found by inverting the parameters,
we denote them $ N', \tilde N'$.

\q If the second (arbitrary ) factors in $N$ and in $  N'$
are taken equal,
so that they cancel out in $N N'^{-1}$, then (with an appropriate
renormalization) we find that the matrix that relates the two products,
namely
$$
N N'^{-1} = R = 1 \oplus (1/q)\pmatrix{ 1 & q-q'\cr0&qq'} \oplus 1,\eqno(5.6)
$$
agrees with (5.3) and satisfies the Yang-Baxter relation.

\ve

\q The formulas that relate the $*$-product to the commutative product
(that one can identify with the ordinary product of functions) depend on the
choice of the arbitrary parameters $A,B$ in the expression for the N-matrices.
In the simplest case, when $A = \tilde A = 1/(q+q'),\, B = \tilde B = 0$,
one has
$$
\eqalign{
&a*b = ab,\, c*a = \,ca,\, c*b = \,cb,\cr
&a*d  = ad + {q-q' \over q+q'}cb.}\eqno(5.7)
$$
It should be noted that, unlike the relations (5.4),
these last formulas do not by themselves define a product on $A$.
They do not, for example, give any information about $a*(bc)$.
There are two ways to complete the picture.
In the first place one can revert to Eq(5.2); this
equation determines a coproduct on $U$
in terms of the standard coproduct $\D_0$, and by duality a product on $A$.
The alternative is to choose a mapping that associates, to each
ordinary polynomial in $a,b,c,d$, a $*$-polynomial. Variations on this theme
were discussed in Section 4D.

\q Formulas of the type (5.2) were introduced by Reshetikhin [5], and used,
with an N-matrix constructed entirely from the Cartan subalgebra of \g,
to discover the generalized quantum groups that are generally known as
twisted versions of the Drinfeld quantum groups.
They can be used to calculate
the further deformations of these twisted quantum groups.
The result can be summarized as an isomorphy of the
cohomology spaces of bialgebras and quantum groups. This completes the
quantization of all the constant r-matrices for the
 simple, complex  Lie algebras.

\bb
\ce{\bf Acknowledgements}
\b
It is a pleasure to acknowledge my indebtedness to Alberto Galindo for
the profitable collaboration on which these lectures are based.
I also thank Miguel Rodrigez, Yat-Lo Lay and
Kim Ming Lau for pleasant and useful discussions and Jerzy Luhierski
for hospitality in Karpacz. I thank the Research
Committee of The Academic Senate, UCLA, for support  and
for a travel grant.

\ve

 \parindent=0.5
cm

\ce{{\bf References}}

\b
\hoffset=0.4truein
\frenchspacing
\item{1.} A.A. Belavin and V.G. Drinfeld,
Sov. Sci. Rev. Math. {\bf 4} (1984), 93-165.
\item{2.} V.G. Drinfeld, {\sl Quantum Groups}, Intern. Congr. Math. Berkeley
1986, 798-820.
\item{3.} M. Jimbo,  Commun. Math. Phys. {\bf 102} (1986), 537-547.
\item{4.} P. Truini and V.S. Varadarajan,
Lett. Math. Phys. {\bf 21} (1991), 287;
P. Bonneau, M. Flato and G. Pinczon, Lett. Math. Phys.
{\bf 25} (1992), 75-84.
\item{5.} A. Sudbery, J. Phys. A{\bf 23} (1990), L697;
A. Schirrmacher, Z. Phys. C{\bf 50} (1991) 321;
N.Y. Reshetikhin, Lett. Math. Phys. {\bf 20} (1990), 331-336.
\item{6.} L. Woronowicz, Publ. Res. Ins. Math. Sci.,
Kyoto Univ., {\bf 23} (1987) 117-181.
\item{7.} N.Y. Reshetikhin, ref. 5.
\item{8.} C. Fronsdal and A. Galindo, {\sl  New Quantum Groups,
Quantum Planes and
   Applications}, preprint UCLA/92/TEP/ and FT/UCM/4/92.

\item{9.} C. Fronsdal and A. Galindo, {\sl The Universal $T$-Matrix},
Proceedings of the 1992 Joint Summer Research Conference
on Conformal Field
Theory, Topological Field Theory and Quantum Groups,
Holyoke, June 1992, to
be published (preprint UCLA/93/ TEP/2);
Lett. Math. Phys. {\bf 27} (1993),
59-71.
\item{10.} C. Fr\o  nsdal, {\sl Universal $T$-matrix
for Twisted Quantum $gl(N)$}, preprint UCLA/93/ TEP/3, to be published in the
Proceedings of the Nato Conference on Quantum Groups, San Antonio, Texas 1993.
\item{11} F. Bayen, M. Flato, C. Fronsdal,
A. Lichnerowicz and D. Sternheimer,
{\sl Deformation Theory and Quantization}, Ann. Phys. {\bf 111}  (1978) 111.
\item{12} A.Connes and M. Flato,
{\sl Closed $*$-products and Cyclic Cohomology}, preprint 1993.
\item{13} P. Bonneau, M. Flato, M. Gerstenhaber and G. Pinczon, {\sl
The hidden group structure of quantum groups: strong duality,
rigidity and preferred deformations}, preprint, Dijon 1993.
\item{14} N.Yu. Reshetikhin and M.A. Semenov-Tian-Shansky, J. Geom.Phys.
{\bf 5} (1988) 533-550.
\item{15} C. Fronsdal, {\sl Some Ideas about Quantization}, Rep. Math Phys.
{\bf 15} (1978) 391.
\item{16} M. Gerstenhaber and S.D. Schack, Contemp.Math.
{\bf 134} (1992) 51-92.
\item{17} C. Fr\o nsdal and A Galindo,
{\sl Deformations of Multiparameter gl(N)}, to appear in Lett.Math.Phys.
\item{18} E.Cremmer and J.-L. Gervais, Commun.Math.Phys. {\bf 134} (1990)
 619-632.
\ve

\bye